\documentclass[11pt]{amsart}
\usepackage{amssymb}
\voffset- .3cm
\numberwithin{equation}{section}

\newcommand{\N}{\mathbb{N}}

\usepackage{amsmath,amssymb}
\usepackage{picins,slashbox,
  graphicx,verbatim,longtable,
  epic,eepic,stmaryrd,rotating}
\usepackage[all]{xy}

\theoremstyle{plain}

\theoremstyle{definition}

\theoremstyle{remark}

\catcode`\@=11
\long\def\@makecaption#1#2{%
    \vskip 10pt
    \setbox\@tempboxa\hbox{
      \small{#1: }\ignorespaces #2}%
    \ifdim \wd\@tempboxa >\captionwidth {%
        \rightskip=\@captionmargin\leftskip=\@captionmargin
        \unhbox\@tempboxa\par}%
      \else
        \hbox to\hsize{\hfil\box\@tempboxa\hfil}%
    \fi}

\newdimen\@captionmargin\@captionmargin=2\parindent
\newdimen\captionwidth\captionwidth=\hsize
\catcode`\@=12

\newlength{\globalparindent}
\begin{document}

\hfill\hfill 
\vskip 1cm

\centerline {\bf  The moduli space of parallelizable 4-manifolds.}

\vskip 3cm

\centerline {\bf Nadya Shirokova.}

\vskip 1cm

\centerline {\bf Abstract}

\vskip.3cm

  In this paper we construct the space of smooth  4-manifolds and find the  homotopy model for  the connected components of the complement to the discriminant. 
  The discriminant of this space is a singular hypersurface and its generic points correspond to manifolds with isolated Morse singularities. 

\vskip .3cm
 
 These spaces can be considered as a natural base for the recent theories studying invariants for families [H], [LL]. In
our applications  we answer the question posed by S. Bauer [B] - we show that the 
theory of Bauer and Furuta can be set up on our configurational space and their invariant is the step-function on chambers.

 We also introduce the definition of the invariant of finite type suggested by the construction of the space and give a simple example of an  invariant
 of order one. 
\vskip 1cm

{\bf 1. Introduction}.

\vskip .3cm

     By studying spaces of objects we implicitly study their diffeomorphism groups and
invariants of finite type. In this paper we construct an infinite-dimensional Euclidean space,
points of which correspond to smooth asymptotically flat parallelizable 4-manifolds.
 Discriminant of
this space, a singular algebraic hypersurface, corresponds  to manifolds with Morse
singularities. Each  4-manifold is represented by  $ |H^1(M^4,Z_2)\oplus H^3(M^4,Z) \oplus \Sigma|$
chambers, where $\Sigma$ is the set of smooth structures on $M^4$. We show that chamber of this space, corresponding to manifold M 
is homotopy equivalent to a $ H^1(M^4,Z_2)\oplus H^3(M^4,Z)$-cover of $BDiff(M^4)$. 
This is a moduli space over which one can consider  family versions of Seiberg-Witten theory.
\vskip .5cm
 In [S1] we constructed several versions of the spaces of 3-manifolds. Manifolds were given as sets of zeros of systems
  of equations (sections) on a trivial normal bundle of a sphere.  Since all 3-manifolds are
parallelizable, i.e. have a trivial tangent bundle, we embedded them into an infinite
dimensional sphere, where they had a trivial normal bundle. Then we extended the system of equations, defining the manifold as its set of zeros
to all of the ambient sphere. 
By the theorem of Mois in dimension 3 smooth, topological and PL categories are equivalent. The case of 4-manifolds is  more complicated since every open 4-manifold carries uncountably many distinct smooth structures.
We restrict ourselves to the class of parallelizable 4-manifolds, all characteristic classes of which vanish.
These manifolds naturally appear, similar to the 3-dimensional case, as the sections of a trivial normal bundle of a sphere.

 All these manifolds are cobordant, have
trivial Euler characteristics  and if they are connected their signature vanishes.  Examples of such manifolds are $T^4$, the quotient of nilpotent Lie group over 
an integer part and more generally $S^1$ - bundles over 3-manifolds.
\vskip .3cm
 
   In his review of the work of Seiberg and Witten, S.Donaldson suggests that by considering r-dimensional family of equations one should get invariants which are cohomology classes in 
   $H^r(BDiff(M))$. 
\vskip .2cm
  The systematic approach to the study of families of invariants was developed  by M. Hutchings, his   "spaces of objects parametrized 
 by a manifold B" can be thought of as chambers of the spaces of 3 and 4-manifolds.  T. Li and A. Liu 
 also    consider family versions of Seiberg-Witten invariants [LL]. We will return to this project in
  the  subsequent paper.
 The trivial loop $\gamma$ in the chamber of the space can be considered as a trivial $M$-bundle over $\gamma$, the loop representing a nontrivial homotopy element of the base will give the monodromy representation of $Diff(M)$.

It would be also interesting to compare our results to the ones of
D.Ruberman [R]. Using 1-parametric version of Yang-Mills equations he studies the topology of the diffeomorphism group. His ``family of objects'' is one-parametric family of metrics (note that the space of metrics
modulo diffeomorphisms  is $BDiff(M)$) , passing  through metric chambers under the action of $Diff$.

\vskip .5cm
  The first application of our approach - is the answer to the question posed by S. Bauer.

  In [B] he gives an interpretation of the monopole equations and the cohomological data of
  the corresponding moduli spaces - Seiberg-Witten invariants in terms of equivariant stable homotopy,
  by considering a monopole map.
   
   He is interested in interpreting the monopole map as a map of infinite-dimensional bundles over the
   configurational space $Conf(M)$ of all choices made: metrics, $spin^c$-connections, harmonic 
   1-forms. Then the monopole map can be understood as an Euler class of the virtual index bundle in
   stable cohomotopy group

   $$\pi_G^0(Conf(M);ind(l))$$

 In this paper we   construct the configurational space $Conf(M)$.
    It can be considered as  the classifying space of G-action (extension of a diffeomorphism group by a gauge group) for parametrized families of 4-manifolds. Therefore
   the image of the monopole class  is the universal parametrized stable cohomotopy invariant which will reflect the information about the diffeomorphism group of 4-manifolds.
  
\vskip .3cm

   Another application of the "space" approach - is the study of invariants of finite type, following the
 original Vassiliev's idea. We give the axiomatics of invariants of finite type and introduce a very simple invariant - the dimension of the moduli
 space of monopoles. This invariant is a constant function on chambers which
 satisfies the condition of being invariant of finite type one.
 \vskip .3cm
 {\bf Acknowledgements}. I want to thank Yasha Eliashberg for important suggestions.

\vskip 1cm

{\bf 2. The space of paralellizable 4-manifolds}.

\vskip .5cm

We consider n smooth function in general position on the trivial normal $n$-bundle to $S^{n+4}$. Generically the common
set of zeros of these functions is a smooth parallelizable 4-manifold, the set of singular immersions
 form a singular hypersurface in the infinite-dimensional space of  functions. The singularities, arising 
 on such discriminant were studied by Goryunov and Mond.

{\bf Definitions.} Denote by $D$ and call the $\bf discriminant$ the subset in $E$ which consists of such tuples $f$ that on the set $\{f=0\}$ there exists a point, s.t. the $f: (S^{n+4},0) \rightarrow (S^4,0)$ is an immersion.
  
  Connected components of $E-D$ we will call $\bf chambers$.
   Part of the discriminant separating two chambers is 
 called a $\bf wall$.

\vskip .2cm
 In this section we modify the previous construction [S2] to get the space of stably parallelizable 4-manifolds.
 
 \vskip .2cm
 
The manifold is parallelizable if it admits the global field of frames, i.e. has
a trivial tangent bundle. In the case of 4-manifolds this condition is equivalent to vanishing of Euler and
the second Stieffel-Whitney class. In particular they all carry $Spin$-structure and the signature of such manifolds is 0.

\vskip .3cm
 By Gromov's h-principle any smooth 4-manifold with all of its smooth structures and metrics can be obtained as a common set of 
zeros of a system of equations in $R^{\N}$ for sufficiently large N. We show that this will hold for any parametrized families of stably paralellizable 4-manifolds, so as the result we get a moduli space of
such manifolds given as sections (systems of equations) in $R^{\N}$:
\vskip .5cm

\vskip .5cm

{\bf Theorem 1.} Any stably parallelizable closed smooth 4-manifold can be obtained as a set of zeros of n functions
on the trivial n-bundle over $S^{n+4}$. Each smooth manifold will be represented by  $|H^1(M,Z_2) \oplus H^3(M,Z)|$ chambers for every smooth structure on this manifold.
\vskip .5cm

{\bf Proof.} First we observe that we need the second Stieffel-Whitney class to vanish in order to have
triviality of the normal bundle and uniqueness of the presentation of our manifold as the set of zeros.
 
We give 4-manifold by equations in a small tubular neighbourhood 
  and then we consider the  obstruction theory problem to extend these equations (sections of the normal bundle to a sphere) to all of  $S^n$ without extra zeros. Namey, we consider the map:
  
  $$ \psi: Framings (M^4) \rightarrow \mathcal Obstructions$$

The homotopy classes of maps from $M$ into the normal bundle, giving the framing,
are found to be
 $$ [M^4:SO(n)]=H^1(M,Z_2)\oplus H^3(M,Z)$$

 To show this we first observe that homotopy groups of $SO(n)$ for large $n$ are known, in particular
  $\pi_1= Z_2$, $\pi_2=\pi_4=0$, $\pi_3=Z$.

 Next we apply Hopf theorem which asserts that for a cell complex one has a one-to-one correspondence between $\pi(M,K(\pi,n))$ and
$H^n(M,\pi)$, homotopy classes of maps into Eilenberg-MacLane spaces and cohomology classes.
\vskip .4cm

 The corresponding obstruction lies in

$$H^i(S^{n+4}, \nu M^4, \pi_{i-1}(S^{n-1}))$$ 

\vskip .3 cm

For dimensional reasons, since all lower homotopy groups are trivial, the only obstruction to extending sections to all of $S^n$ will lie in

$$H^{n+4}(S^{n+4}, \nu M^4, \pi_{n+3}(S^{n-1}))=\pi_{n+3}(S^{n-1})$$ 

$$\pi_{n+3}(S^{n-1})=0$$
\vskip .5cm

 Thus the obstruction map becomes
 
 $$\psi :  H^1(M^4,Z_2)\oplus H^3(M^4,Z) \rightarrow 0$$
 
 This implies, that the section (system of equations) given in  tubular neighborhood can be extended to all of $S^n$ for any
 framing we start with.
 
 In particular we will be getting all of $M^4=S^1 \times V^3$ will all smooth structures. This is an important class of 4-manifolds
posessing many useful properties.

\vskip .4cm

{\bf Proposition}. If $M^4=S^1 \times V^3$, M is spin, $b_2(M)$ is even and $b_+(M)=b_-(M)=b_1(V)$.

\vskip .2cm

\vskip .2cm

 When $M^4=S^1 \times V^3$ the set of framings has the most interesting (from the point of view of Seiberg-Witten theory) form:

$$Z\oplus Spin^c(V)$$
i.e. each such manifold will be represented with all of its $Spin^c$-structures.
\vskip .5cm

We can also consider punctured asymptotically flat 4-manifolds. By the theorem of Quinn, any  punctured 4-manifold posesses a smooth structure (and even an infinite number on them). This will
eliminate cases when we don't know if a 4-manifold has any smooth structures at all. (Such as $S^2 \times S^2$). But now every manifold will come up with uncountably many smooth structures.

  However by the theorem  of Vidussi  which shows that the
manifolds diffeomorphic outside a point have the same Seiberg-Witten 
invariants (so one cannot use them to detect eventual inequivalent smooth 
structures), the parametric families of Seiberg-Witten invariants of such manifolds will be the same.
 
\vskip .5cm
 We make the following observation:
 \vskip .3 cm
{\bf  Lemma1. } Stably parallelizable  open 4-manifolds are parallelizable.
\vskip .3 cm
{\bf Proof.} The sum of the tangent bundle and the copy of $\mathbf R$ is trivial. Since the base
doesn't have the highest homology group, Euler class is trivial and the tangent bundle itself can be trivialized.

 \vskip .2cm

If will work with ``asymptotically flat'' 4-manifolds, i.e. such that outside the 
ball $B_R$ of some
large radius R they will be given as the set of common zeros of the system of linear equations (e.g. $f_i(x_1,...x_{n+4})=x_i$ for $i=1,...n$.),  the analogue of Theorem 1 will hold as well
 and the proof  is the same:

\vskip .5cm

{\bf Theorem 1'.} Any smooth, parallelizable, asymptotically flat  4-manifold can be obtained as a set of zeros of n functions
on the trivial (n+4)-bundle over $S^n$. Each smooth manifold will be represented by  $|H^1(M,Z_2) \oplus H^3(M,Z)|$ chambers for every smooth structure on this manifold.

\vskip .5cm
{\bf 3. The topology of the chambers of the spaces of manifolds.}

\vskip .5cm
 In this section we set n to be infinite, i.e. we embed parallelizable 4-manifolds into infinite-dimensional sphere and represent  them by  systems of equations.
 
 The following theorem is supposed to be known to specialists, however since I was not able to find
 its proof in the literature, I  sketch it here. The same proof will work for manifolds of any dimension,
 embedded into $S^{\infty}$.

\vskip .2cm

{\bf Theorem 2.} Chambers of the space of 4-manifolds corresponding to manifold
of homemorphism type M are homotopy equivalent to a $|H^1(M,Z_2)\oplus H^3(M,Z)|$-cover of
$ BDiff(M)$.

\vskip.3cm

  Each chamber of the space consists of all possible (nonparametrized embeddings of $M^4$ into $S^n$. We consider a ``representative'' for the given homoeomorphism class
( perhaps with an additional structure such as framing or homology class), a manifold $M_0$,  given by a 
simple system of equations. Over the chamber we consider a bundle with fiber Diff(M), all
 diffeomorphisms from $M_0$ to the manifold, corresponding to the point of the chamber. This
is a principle bundle with the structure group Diff(M). If we show that the total space of the bundle is contractible
 and the action of
Diff(M) is free, then by the definition of the classifying space, our chamber is homotopy equivalent to BDiff(M).

 Theorem will hold if we embed a 4-manifold into an
infinite-dimensional space and obtain it as set of zeros on n functions
of (n+4) variables, where n approaches infinity.

\vskip .2cm

 {\bf Remark 1}. Note that in the case of 3-manifolds A. Hatcher and D. McCullough answered
 the question posed by M. Kontsevich, regarding the finiteness
 of the homotopy type of the classifying space of the group of
 diffeomorphisms [HM]: 
\vskip .2cm
{\bf Theorem [HM]}  Let M be an irreducible compact connected orientable
3-manifold with nonempty boundary. Then $BDiff(M, rel \partial)$ has the homotopy type of an aspherical CW-complex.

 It would be very interesting to see what could be the analogue of the Hatcher's theorem in the case of
4-manifolds.

\vskip .5cm

We are looking for a simple topological space to which the chamber $C_M$ of the space of
 4-manifolds is homotopy equivalent.

 In our proof  we use the analogy between the construction of the Stieffel
manifold $St(m,\infty)$, the universal bundle for the Grassmanian $Gr(m,\infty)$
 and the space of 
embeddings of a 4-manifold into an infinite-dimensional Euclidean space.
Let us denote by $E$ the space of all embeddings of a ``standard'' manifold $M_0$ given by a fixed system of equations into
an infinite-dimensional space $E=Emb(M_0^4, R^{\infty)}$; then there is a natural action of the group
$Diff(M_0)$. Let's consider the quotient $E/Diff = N$. E is a principal Diff-bundle over N. 

Let's also denote by $Emb(L_0,L)$ the set of embedding of the standard m-plane into
a given one.
\vskip .5cm

 In the following table we  show the analogy between Grassman and Stiffel manifolds
and their nonlinear analogues - spaces of $M^4$ in $R^{\infty}$, parametrized by equations.

\vskip .3cm
 
  \begin{center}
     \begin {tabular}{|l|r|}
        \hline
      $St(m,\infty)$ & E=all parameterized embeddings of $M^4$ into $R^{\infty}$ \\ \hline
      Grassmanian=all linear m-subspaces of $R^{\infty}$ & Chamber=L, embeddings of $M^4$ into $R^{\infty}$ \\ \hline
      St=principal O(m)-bundle over Gr  &  E=principal Diff(M)-bundle over L \\ \hline
      Fiber=set of frames in a given plane & Fiber=set of all diffeo's from a given $M_0$ into M \\ \hline
      Frame=emb. of a standard plane into a given one & Diffeo from $M_0$ to M  - ``change of variables'' \\ \hline
      Total space=Stiffel manifold & Total space=E, parametrized embeddings \\ \hline
      $St(m, \infty)$ is contractible & E is contractible (Lemma) \\\hline
      \end{tabular}
  \end{center}

\vskip 1cm
{\bf Remark 2}. The above analogy suggests that one can think of the space of 4-manifolds as a nonlinear analogue of
 the Grassmanian of 4-planes in $R^{\infty}$.

\vskip .3cm
{\bf Lemma 2.}  The space E is contractible, the action of the group of
diffeomorphisms on M is free.

\vskip .2cm

{\bf Proof}. We know by Whitney's theorem any two embeddings of an n-dimensional 
manifold into a $(2n+2)$-dimensional Eucledean space are isotopic. This means that the 
space E is at least connected. To see that it is contractible we show that
all of its homotopy groups are trivial. This can be proved by showing that any two isotopies
 can be connected by an isotopy, then two isotopies between isotopies can also be connected by an isotopy.
If the ambient space is infinite-dimensional, we can kill all homotopy groups by iterating this process.

 An isotopy between
two distinct embeddings of the same manifold M is a map from $M \times I$
into $R^{2n+2}$, the image of which is not necessarily embedded. However, 
it always can be embedded into $R^{2n+4}$, i.e. by increasing the dimension 
of the ambient manifold we can assume this. Now these two isotopies are
isotopic by Whitney. Continuing this process and adding the dimensions to
the ambient manifold, we prove that the space of
the embeddings of a given 4-manifold M into an infinite-dimensional Euclidean
space is contractible. 

 The action of Diff(M) on E is free by the construction.

\vskip .2cm

 The chambers of the spaces of smooth 4-manifolds that we constructed in section 1 correspond not
just to all possible embeddings of M into $S^ \infty$, but to M with  framing . So we have to consider a fibration of the
group Diff(M) by subgroups, preserving a given structure.

\vskip .2cm 

{\bf Lemma 2}. Consider the 
subgroups $Diff_{fr_i}(M)$ in $Diff(M)$, preserving a given framing $fr_i$. For any such 
subgroup we get a fibration:

 $$BDiff(M) {\buildrel {Diff(M)/Diff_{{fr}_i}(M)} \over {\longrightarrow}} BDiff_{{fr}_i}$$

with the fiber -  $|H^1(M,Z_2)\oplus H^3(M,Z)|$.

\vskip .2cm

{\bf Proof}. The set of normal framings on a 4-manifold by the definition is [M:SO(n)],
 which represent homotopy classes of trivializations of the normal bundle and which is
 found to be $H^1(M,Z_2)\oplus H^3(M,Z)$, or  in the case of $M^4=M^3\times S^1$ equals $Z^2\oplus Spin\oplus Spin^c$. 
Thus factoring the group of diffeomorphisms of
a 4-manifold by the subgroup, preserving one framing, we get this set as the quotient.
Here $Diff(M)$ is the group of diffeomorphisms of any given smooth structure on $M$.

\vskip .4cm
{\bf Remark 3}.  It is useful to think of $BDiff(M)$ as a (contractible) space of metrics modulo diffeomorphisms.

\vskip .5cm
{\bf 4. The embedded 5-cobordisms.}

\vskip .5cm

 Here we show that the 5-cobordisms, connecting 4-manifolds are embedded into $(S^{\infty} \times I)$
 and given by the equations.
 \vskip .5cm
 {\bf Theorem 3}. If $(X_1^4, f_1)$ is spin cobordant to $(X_2^4, f_2)$ in $(S^{\infty} \times I)$, then
 $\psi (X_1^4, f_1) = \psi (X_2^4, f_2)$, where $f_i$ is a framing of $X_i$.
 \vskip .3cm
 {\bf Proof}. Let $W^5$ be spin cobordism  $$\partial W^5 = (X_1^4, f_1) \cup (X_2^4, f_2)$$
 
 the obstruction to extending section from a tubular neighborhood of $W^5$ lies in
 
 $$H^{n+5}(S^{n+4}, \nu W^5, \pi_{n+4}(S^{n-1}))=\pi_{n+4}(S^{n-1})$$ 
we observe that
$$\pi_{n+5}(S^n)=0$$

Thus we get a commutative diagram:

 \vskip 1cm

\begin{equation*}
\xymatrix@C+0.5cm{f \in [W^5 : SO(n)]
\ar[d]
\ar[r] & f_1 \in [X_1 : SO(n)]  \ar[d]^{\psi_1} \\
f_2 \in [X_2 : SO(n)]  \ar[r]^{\psi_2} & \pi_{n+4}(S^{n-1})}
\end{equation*}
 \vskip 1cm
 
 Once $W^5$ exists,  $\psi_1$ and $\psi_2$ map to the same element of $\pi_{n+5}(S^n)$. This
 shows the naturality of the primary obstruction.
\vskip .5cm

{\bf Conjecture. }  The space of 4-manifolds is the classifying space of the category of embedded 5-cobordisms.
\vskip .5cm

\vskip 1cm
{\bf 5. The hyperkahler structure.}
\vskip .7cm

 In this paragraph we show that all  parallelizable open 4-manifolds carry the hyperkahler structure. 
 \vskip .3cm
 A hyperkahler manifold is a Riemannian manifold of dimension $4k$ with holonomy group contained in $Sp(k)$ (the group of quaternionic-linear unitary endomorphisms of an n-dimensional quaternionic Hermitian space). 
 
 Non-compact hyperkahler 4-manifolds which are asymptotic to $H/G$, where H denotes the quaternions and $G$ is a finite subgroup of $Sp(1)$, are known as Asymptotically Locally Euclidean, or ALE, spaces. These space are studied in physics and are called the  gravitational instantons.

 \vskip .5cm
{\bf Theorem 3}. Every open parallelizable 4-manifold carries the hyperkahler structure.
\vskip .3cm
{\bf Proof}.  It is a well-known fact that any even dimensional vector space admits a linear complex structure. Therefore an even dimensional manifold always admits a $(1,1)$ rank tensor pointwise (which is just a linear transformation on each tangent space) such that $J^2 = 1$ at each point. 

 When this local tensor can be defined globally, the pointwise linear complex structure yield an almost complex structure, which is then uniquely determined. The existence of an almost complex structure on a manifold M is equivalent to a reduction of the structure group of the tangent bundle from GL(2n, R) to GL(n, C) and it is possible because of the triviality of the characteristic classes. Now we want to integrate this structure to a complex one.

 Recall, that according to the theorem of Landweber [L],  every open even dimensional manifold
which has trivial homology above the middle dimension, carries the complex structure.
For open 4-manifold, the highest homology group is trivial, but generically $H^3(M) \neq 0$. If we
 extend the complex structure to the 3-handles, then open $M^4$ will be complex. It is possible since the
 manifold is parallelizable, i.e. we can choose a global parallelization a uniform choice of frames and the  Nijenhuis tensor  vanishes. So all open parallelizable 4-manifolds carry complex structure.
 \vskip .3cm
 Next we immerse open complex $M^4$ into ${\mathbf C}^2$ , s.t. the Jacobian of the map is everywhere
 nondegenerate which is always possible to arrange by the Smale-Hirsch theorem.
 
 It follows, that all the structures one has on ${\mathbf C}^2$ are pulled back on $M^4$  under this immersion.
 
 In particular we get that every open parallelizable 4-manifold carries kahler and  the hyperkahler
 structures. We are mostly interested in those hyperkahler metrics, which are complete.

 It would be very interesting to understand how these " hyperkahler subchambers" located in the
 space of all metrics (chamber of the space of manifolds), see the Remark 3.

\newpage
{\bf 6. Refined SW-invariants.} 
\vskip .7cm
 The theory which naturally fits into our space  is the "refined SW-invariants"
 defined in purely topological terms via stable cohomotopy. The purpose of this paper is to construct  the space of 4-manifolds, so here we just explain the relation to the Bauer-Furuta theory. The detailed
discussion will appear in the subsequent paper.

\vskip .5cm

 Stable cohomotopy invariants were introduced by Bauer and Furuta [BF] and are defined via
 a version of Pontryagin-Thom construction by considering the monopole map instead of monopole
 equations.

The monopole map is a $U(1)$-equivariant map between affine Hilbert spaces. The moduli space of monopoles is obtained as a quotient of the 
zero-set of the monopole map by the action of the group $U(1)$ of complex numbers of unit length.
 In the finite-dimensional case a proper map between two spaces can be extended to a one-point compactification and then by the
Pontrijagin-Thom construction will define an element in the stable homotopy groups of spheres.
By the theorem of Schwarz [Sc] this construction can be carried out in the infinite-dimensional setting.

$$\mu  : Conn \oplus (\Gamma (S^{+})  \oplus \Omega^1(M)/G \rightarrow Conn \oplus (\Gamma(S^+) \oplus \Omega^+(M) \oplus H^1(M,R) \oplus \Omega^0(M)/G)$$

$$(A,\psi,a) \rightarrow(A,D_{A+a} \psi, F^+_A-\sigma(\psi), a_{harm},d^*a)$$

As a map over $Conn$, space of $spin^c$-connections, monopole map is equivariant with respect to the action of
$G=map(M,U(1))$. If A-fixed connection, then $Conn$ is invariant under the free action of gauge group
with the quotient space isomorphic to

$$Pic^0(M)=H^1(M,R)/H^1(M,Z)$$

\vskip .5 cm

{\bf Theorem [B-F]}.  The monopole map $\mu$ for an oriented 4-manifold M with $spin^c$-structure s defines an element in the equivariant stable
cohomotopy group

$${\pi^b}_{U(1)}(Pic^0(M);ind(D)$$

which is independent of the chosen Riemannian metric. For $b \ge dim(Pic^0(M)) + 1$, a homology
orientation determines a homomorphism of this stable cohomotopy group to $Z$, which maps $[\mu]$
to the integer valued Seiberg-Witten invariant.

\vskip .5 cm
 Here $ Pic^0(M)$ denotes a Picard torus. The Dirac operator associated to the chosen $spin^c$-connection defines a virtual complex bundle $ind(D)$ over the Picard torus, and $b=b_+(M)$. By
 the homology orientation we understand the orientation of $H^1(M,R) \oplus {H_+}^2(M,R)$ and the
 index of the Dirac operator is a complex vector space of complex dimension
 $$d=1/8(c_1(s) - \sigma (M))$$
 
 where $c_1(s)$ is the Chern class of the spinor bundle.
 \vskip .5cm
  The next theorem/observation describes the space  $Conf(M)$ , which consists
 of all choices made: metrics, $Spin$-connections, harmonic 
   1-forms. The monopole map can be then understood as an Euler class of the virtual index bundle in
   stable cohomotopy group

   $$\pi_G^0(Conf(M);ind(l))$$
   
   where $G$ is the extension of the subgroup of $Diff_{fr}(M)$ (preserving the homological orientation)
   by the gauge group $U(1)$.
   
   The explicit calculation of this invariant for parametric families will lead to the information about the topology of the diffeomorphism group of a four-manifold.
   
   \vskip .5cm

    We can organize all "data" coming from monopole equation and from the space of 4-manifolds
    in the following theorem:
   \vskip .3cm
   
 {\bf Theorem/Definition 4.}  The monopole map defines an element in an equivariant stable cohomotopy
 group
  $$\pi_G^0(Conf(M_{fr});ind(l))$$ where G is an extention of the group of diffeomorphisms preserving the homological orientation
  by the gauge group. For a given stably parallelizable 4-manifold M the union of spaces $Conf(M_{fr})$ over all framings is the space of metrics
  on M times the set of spin-structures, holomorphic 1-forms and the gauge equivalence classes of
  connections in $U(1)$-bundle.
   \vskip .3cm
   
   {\bf Proof}. Over each parallelizable 4-manifold we consider principal $U(1)$-bundle.  Denote as $\mathcal B$ the total space of this bundle over the chamber.
    As before, consider the space of connections in $U(1)$-bundle up to the gauge equivalence.

     Denote by $C_{M_{fr}}$ the  chamber of the space of stably parallelizable 4-manifolds, corresponding to manifold $M_{fr}$, constructed in paragraph 2 and the gauge equivalence classes of connections in the total space of the $U(1)$-bundle  $\mathcal B$, as
     $Conf(M_{fr})$. 
 By the Theorem 2 the chamber, corresponding to manifold M is homotopy equivalent to a  $|H^1(M,Z_2)\oplus H^3(M,Z)|$-cover of $BDiff(M)$, i.e. this is a classifying space of the subgroup of the group
 of diffeomorphisms preserving cohomology class from  $H^1(M,Z_2)\oplus H^3(M,Z)$. 
 We also have shown in paragraph 3 that  $BDiff(M)$ can be considered as a space of metrics up to diffeomorphism. So the union of all chambers, corresponding to M is homotopy equivalent to the 
 space of metrics.
 
  Next we will fit  into our picture the other parameters involved in the definition of the moduli space of monopoles:
     \vskip .3cm
  
  $H^1(M,Z)$ can be identified with the set of harmonic 1-forms, i.e. there is one harmonic representative 
in each cohomology class (and it is dual to  $H^3(M,Z)$).

   \vskip .3cm
 $H^1(X,Z_2)$ - is the set of $spin$-structures, or $Z_2^{b_1(M)}$. Note that given a $spin$-structure 
 on a manifold, one can find the corresponding $spin^c$-structure.
   \vskip .3cm

 The monopole map is defined as follows:

$$(g_{i,j},A,a,\phi) \rightarrow(g_{i,j},A,a,D_{A+a} \phi, F^+_A-\sigma(\phi), a_{harm},d^*a)$$

Then as before we have the parametrized monopole map:
    
   $$\mu  : Conf(M_{fr}) \oplus (\Gamma (S^{+})  \oplus \Omega^1(M)/G \rightarrow Conf(M_{fr}) \oplus (\Gamma(S^+) \oplus \Omega^+(M) \oplus H^1(M,R) \oplus \Omega^0(M)/G$$

\vskip 1cm

{\bf 7. Example of an  invariant of order one.}

\vskip .5cm
 The construction of the space of knots led V. Vassiliev to the axiomatics of invariants of finite type.
 He introduced a spectral sequence which calculated the cohomology of the complement of the
 discriminant of the space of knots. The filtration in this spectral sequence produced the notion
 of the invariant of order n.
 
  It is tempting to try this approach in the case of higher-dimensional manifolds. Then the topology
  of the space of objects should dictate the natural definition of the invariant.

   \vskip .3cm
   {\bf Definition.} An invariant of an asymptotically flat parallelizable 4-manifold is of {\bf finite type n} if for
   any selfintersection of the discriminant of codimension $n$ its alternated sum over $2^n$
   chambers adjacent to this selfintersection is zero.
   
     \vskip .5cm
 We constract first simple (scalar) invariant of order one. 
 
 Consider the moduli space of monopoles of a closed parallelizable 4-manifold.
 Then its dimension is given by the formula:
 
 $$dim  \mathcal{ M} =1/4(c_1^2(M) - 2 \chi (M) - 3 \sigma (M))$$
 
 Of course, this invariant doesn't distinguish smooth structures and is the same on chambers, 
 corresponding to different diffeomorphism types of the manifold of a given homotopy type M.
 \vskip .5cm
 
 {\bf Theorem 5}. $dim \mathcal{M}(mod 2)$  is an invariant of a stably parallelizable 4-manifold $M^4$ of order one.

 \vskip .5cm
 
 {\bf Proof.}
 Recall that the Chern class and the Euler characteristics of a spin manifold are divisable by two (all  parallelizable
 manifolds are spin). Signature  of a stably parallelizable manifold is zero. Thus our invariant is just
 an integer, depending on the Chern class and the Euler characteristics. We want to show that given
 any four chambers, adjacent to the selfintersection of the discriminant of codimention two, the
 alternated sum of the invariants, corresponding to these chambers is zero.
 
  Indeed, every wall of the discriminant corresponds either to a surgery  or to the change of a spin-structure, both of which changes $\chi (M)$ by
  two (and an invariant by one), or to the change of smooth structure, which doesn't change the invariant Given selfintersection of the
  discriminant of codimension two, there are four chambers, adjacent to this selfintersection. Our
  invariant changes at most by one by passing from one chamber to another, i.e. the alternated sum
  of invariants modulo 2 is zero.
 
\newpage
 
{\bf 8. Further directions.}

\vskip .7cm

 1. We would like to consider a parametric versions of Seiberg-Witten equations expecting to construct  "secondary Seiberg-Witten invariants".
 Over each manifold in a chamber, corresponding to a given  4-manifold we consider a $U(1)$-bundle and connection
and fix a family of fiberwise  $Spin^c$-structures. Then the space of fiberwise monopoles gives rise to the
family Seiberg-Witten invariants.
 Since the homotopy type of the chamber of our space is $BDiff(M^4)$ by considering r-dimensional families of equations
we should encount, in particular $H^r(BDiff(M))$. We see that in the case of $M^4=M^3\times S^1$ we can get a moduli space of dimension -1 and apply
 Ruberman's technique to get invariants of the diffeomorphism group. Since all parallelizable manifolds
 are not simplyconnected, this construction should provide an insight into the topology of these manifolds and their diffeomorphism groups.
\vskip .3cm
2. Our next goal is to consider the local system of Gukov-Witten [GW] invariants on our space, extending the program, outlined in [S2] to 4-manifolds.

 Recently S.Gukov and E.Witten  introduced a categorification of Vafa invariant [GW]. This theory
 fits into our program and using the local system of their complexes on the space of parallelizable
 4-manifolds we can construct the Gukov-Witten functor. This construction should be useful for the
 study of the diffeomorphism group of 4-manifold. Even the case of the 4-sphere is still open.

Once such sheaf is constructed, the homology of the 
group of diffeomorphisms will be included into its homology.
All these developments suggest that the corresponding sheaf should contain the information about
diffeomorphisms of 4-manifolds and a variant of Vassiliev theory.

\vskip .3cm
3. In order to define Vassiliev-type invariants for 4-manifolds it would be very interesting to extend this local system to the singular locus. This construction
should provide invariants of 4-manifold with Morse singularity or a "Vassiliev derivative" for the
refined Seiberg-Witten invariant. 

\vskip .3cm
4. We conjecture that the space which we constructed is the classifying space of the category of embedded 5-cobordisms.

\vskip .3cm
5. It would be very interesting to relate our approach to the one of  P. Kronheimer  and  C. Manolescu
[M], [KM].
\vskip .3cm
6. We would like to understand which of the hyperkahler metrics that we constructed are complete and how they are located in $BDiff(M)$.

\vskip 1cm

{\bf 9. Bibliography.}
\vskip 1cm

[B] S. Bauer, Refined Seiberg-Witten invariants, GT/0312523
\vskip .2cm
 [BF]  S. Bauer, M. Furuta, A stable cohomotopy refinement of Seiberg-Witten invariants: I. math.DG/0204340
 
 \vskip .2cm
[GW]  S. Gukov, E. Witten, Notes on Gauge Theory and Categorification, preprint 2006.
\vskip .2cm
[HM] A. Hatcher, D. McCullough,Finiteness of Classifying Spaces of Relative Diffeomorphism Groups of
 3-manifolds,Geom.Top., 1 (1997).

\vskip .2cm
[Hu] M. Hutchings,Floer homology of families 1, preprint SG/0308115.
\vskip .2cm
[K] P. Kronheimer, J.Diff. Geom., 29, 1989, no.3, 665-683. 
\vskip .2cm
[KM]  P. Kronheimer, C. Manolescu,  Periodic Floer pro-spectra from the Seiberg-Witten equations, arXiv:math/0203243. 
\vskip .2cm
[L] P.Landweber, Complex structures on open manifolds, Topology, Vol.13, p.69-75, 1974.
\vskip .2cm
 [LL] T.J.Li, A.K.Liu, Family Seiber-Witten invariants and wall crossing formulas, preprint GT/0107211

\vskip .2cm
[M]  C. Manolescu, Seiberg-Witten-Floer stable homotopy type of three-manifolds with $b_1$=0, Geom. Topol. 7(2003) 889-932.

\vskip .2cm
[R] D. Ruberman,A polynomial invariant of diffeomorphisms of 4-manifold, geometry and Topology, v.2, Proceeding of the Kirbyfest,
p.473-488, 1999.
\vskip .2cm
[S] A. Schwartz,Dokl.Acad.Nauk USSR, 154, pp.61-63.
\vskip .2cm
[S1]  N. Shirokova, The space of 3-manifolds, Comptes Rendus Acad.Sci., t.331, p. 131-136, 2000.
\vskip .2cm
[S2] N. Shirokova, On the classification of Floer-type theories., arXiv:0704.1330.

\end {document}